\documentclass[12pt]{article}
\newtheorem{lemma}{Lemma}
\newtheorem{theorem}[lemma]{Theorem}

\def\deg{{\rm deg}}
\def\qed{\hspace*{\fill}\vrule height6pt width4pt depth0pt\medskip}

\begin{document}
\title{{\scshape\normalsize Mathematics Division, National Center for
         Theoretical Sciences at Taipei}\\
         {\scshape\large NCTS/TPE-Math Technical Report 2004-005} \\~\\
         {\bf A note on equitable colorings of forests}\thanks{This is
                a note for a talk of the Seminar in Discrete Mathematics
                (NCTS/TPE and Academia Sinica) in March 26, 2004.}}
\author{Gerard J. Chang\thanks{Department of Mathematics,
                National Taiwan University, Taipei 10617, Taiwan.
                Email: gjchang @math.ntu.edu.tw.
                Support in part by the National Science Council
                under grant NSC92-2115-M002-015.
                Member of Mathematics Division, National Center for
                Theoretical Sciences at Taipei.}}
\date{April 13, 2004}
\maketitle

\begin{abstract}
This note gives a short proof on characterizations of a forest to
be equitably $k$-colorable.
\end{abstract}

\section{Introduction}

In a graph $G=(V,E)$, a {\it stable set} (or {\it independent
set}) is a pairwise non-adjacent vertex subset of $V$.  The {\it
stability number} (or {\it independence number}) $\alpha(G)$ of
$G$ is the maximum size of a stable set in $G$.
An {\it equitable $k$-coloring} of $G=(V,E)$ is a partition of $V$
into $k$ pairwise disjoint stable sets $C_1$, $C_2$, $\ldots$,
$C_k$ such that $||C_i|-|C_j|| \le 1$ for all $i$ and $j$. The
{\it equitable chromatic number} $\chi_=(G)$ of $G$ is the minimum
number $k$ for which $G$ has an equitable $k$-coloring.

The notion of equitable colorability was introduced by Meyer
\cite{m}, who also conjectured a statement stronger than Brooks'
theorem that $\chi_=(G) \le \Delta(G)$ for any connected graph $G$
other than a complete graph or an odd cycle, where $\Delta(G)$ is
the maximum degree of a vertex in $G$.  Hajn\'al and Szemer\'edi
\cite{hz} gave a deep result that any graph $G$ is equitably
$k$-colorable for $k > \Delta(G)$.  This topic is then studied for
many researchers. Lih \cite{l} gave a survey on this line.

The main concern of this note is on the equitable colorability of
trees. Meyer in his paper \cite{m} also showed that a tree $T$ is
equitably $(\lceil\frac{\Delta(T)}{2}\rceil+1)$-colorable.
However, this proof was faulty.  It was reported by Guy \cite{g}
that Eggleton remedied the defects. He could prove that a tree $T$
is equitably $k$-colorable if $k \ge
\lceil\frac{\Delta(T)}{2}\rceil+1$. Meyer's results on trees was
greatly improved by Bollob\'as and Guy \cite{bg} as follows.

\begin{theorem} {\bf (Bollob\'as and Guy \cite{bg})}
A tree $T$ of order $n$ is equitably $3$-colorable if $n \ge
3\Delta(T)-8$ or $n=3\Delta(T)-10$.
\end{theorem}

Using this result as the induction basis, Chen and Lih \cite{cl}
gave a complete characterization for a tree to be equitably
$k$-colorable.  Their results are in two parts.  Notice that as a
tree is a connected bipartite graph, its vertex set has a
bipartition.

\begin{theorem}{\bf (Chen and Lih \cite{cl})}
Suppose $T$ is a tree of order $n$, and $(A,B)$ is a bipartition
of $T$. For $||A|-|B|| \le 1$, the tree $T$ is equitably
$k$-colorable if and only if $k \ge 2$.
\end{theorem}

To see their second result, we need another notion. Suppose $x$ is
a vertex in a graph $G=(V,E)$.  An {\it $x$-stable set} in $G$ is
a stable set which contains $x$. The {\it $x$-stability number}
$\alpha_x(G)$ of the graph $G$ is the maximum size of a $x$-stable
set in $G$.  We use $\alpha_x$ for $\alpha_x(G)$ when there is no
ambiguity on the graph $G$.

Suppose $x$ is a vertex in a graph $G=(V,E)$ of order $n$.
Partition $V$ into $k=\chi_=(G)$ stable sets $C_1$, $C_2$,
$\ldots$, $C_k$ such that $||C_i|-|C_j|| \le 1$ for all $i$ and
$j$. Suppose $x\in C_i$. Then $|C_i| \le \alpha_x$ and $|C_j| \le
\alpha_x+1$ for all $j\ne i$.  Consequently,
$$
     n = \sum_{i=1}^k |C_i| \le \alpha_x + (k-1) (\alpha_x+1)
     = \chi_=(G) (\alpha_x+1) - 1
$$
and so $\chi_=(G) \ge \frac{n+1}{\alpha_x+1}$, which gives (see
\cite{cl})
\begin{eqnarray}  \label{lower bound}
  \textstyle{    \chi_=(G) \ge \max\limits_{x\in V}
                 \lceil \frac{n+1}{\alpha_x+1} \rceil.}
\end{eqnarray}

\begin{theorem}  {\bf (Chen and Lih \cite{cl})}     \label{cl}
Suppose $T$ is a tree of order $n\ge 2$, and $(A,B)$ is a
bipartition of $T$. For $||A|-|B||\ge 2$, the tree $T$ is
equitably $k$-colorable if and only if $k \ge \max\{3,
\lceil\frac{n+1}{\alpha_v+1}\rceil\}$, where $v$ is an arbitrary
vertex of degree $\Delta(T)$.
\end{theorem}

Notice that when $||A|-|B|| \le 1$, it is the case that $\alpha_x
\ge \frac{n-1}{2}$ and so $\frac{n+1}{\alpha_x+1} \le 2$ for any
vertex $x$. On the other hand, even when $||A|-|B|| \ge 2$, it is
still possible that $\frac{n+1}{\alpha_x+1} \le 2$ for all
vertices $x$. An easy example is the tree obtained from a $3$-path
by adding $\ell \ge 3$ leaves joining to each vertex of the
$3$-path. This shows that the $3$ in the lower bound of Theorem
\ref{cl} can not be dropped.

An unpublished manuscript by Miyata, Tokunaga and Kaneko
\cite{mtk} gave another characterization of equitable colorability
of trees.  While the proof is long, it is without using other
results.

\begin{theorem}  {\bf (Miyata, Tokunaga and Kaneko \cite{mtk})}
\label{ntk}
 Suppose $T=(V,E)$ is a tree of order $n$ and $k\ge 3$ is an integer.
 Then $T$ is equitably $k$-colorable if and only if
$\alpha_x \ge \lfloor \frac{n}{k}\rfloor$ for any vertex $x$ or
equivalently $k \ge \max_{x\in V}\lceil\frac{n+1}{\alpha_x
+1}\rceil$.
\end{theorem}

Notice that the equivalence follows from that
$$
    {\textstyle k \ge \lceil \frac{n+1} {\alpha_x + 1}\rceil \iff
     k \ge \frac{n+1}{\alpha_x+1} \iff
     \alpha_x \ge \frac{n+1}{k} -1 =\frac{n-k+1}{k}
     \iff \alpha_x \ge \lfloor\frac{n}{k} \rfloor.
    }
$$

 The purpose of this note is to clarify the relation
between Theorems \ref{cl} and \ref{ntk}.  We also give a short
proof of the result by combining all techniques in
\cite{bg,cl,mtk} together. We present the proof in terms of
forests as it is the same as that for trees.

\section{Equitable coloring on forests}

We first clarify the relation between taking maximum over all
vertices in Theorem \ref{ntk} and using only one vertex in Theorem
\ref{cl}.

In a graph $G$, the {\it neighborhood} $N(v)$ of a vertex $v$ is
the set of all vertices adjacent to $v$, and the {\it closed
neighborhood} $N[v]$ is $\{v\} \cup N(v)$. For a subset $S$ of
vertices, the neighborhood  $N(S)$ of $S$ is $\cup_{v\in S} N(v)$.

\begin{lemma} \label{major vertex}
Suppose $v$ is a vertex in a forest $F=(V,E)$ of order $n$.  If
$\lceil \frac{n+1}{\alpha_v + 1}\rceil > 3$, then $v$ is the only
vertex of degree $\Delta(F)$. Consequently, if $\max\{3,
\max_{x\in V} \lceil\frac{n+1}{\alpha_x+1}\rceil \} > 3$, then the
maximum is attained by the unique vertex of degree $\Delta(F)$.
\end{lemma}
\noindent {\bf Proof.} Notice that $\lceil \frac{n+1}{\alpha_v +
1}\rceil > 3$ implies $\frac{n}{\alpha_v+1} \ge 3$ or $n \ge
3\alpha_v+3$.
Suppose $v$ is of degree $d$.  First, $\alpha_v = 1 +
\alpha(F-N[v])$. Notice that the stability number of any bipartite
graph is at least the half of its order as the larger part in a
bipartition is a stable set. It is then the case that $2 \alpha_v
= 2 + 2 \alpha(F-N[v]) \ge 2 +
     n-1-d \ge 2 + 3\alpha_v+3-1-d$
and so $\deg(v) = d \ge \alpha_v +4$.
On the other hand, suppose $x$ is a vertex other than $v$. Then
all of its neighbors, except possibly one, form a stable set in
$F-N[v]$ since $F$ has no cycles.  Hence, $\alpha(F-N[v]) \ge
\deg(x)-1$ and so $\alpha_v = 1 + \alpha(F - N[v]) \ge \deg(x)$,
which in turn implies $\deg(v) > \deg(x)$. \qed

Lemma \ref{major vertex} implies that the conditions in Theorems
\ref{cl} and \ref{ntk} are in fact the same. Having this in mind,
we are ready to re-prove the main assertion.

\begin{theorem}
 Suppose $F$ is a forest of order $n$ and $k \ge 3$ is an integer.
 Then $F$ is equitably $k$-colorable if and only if
$\alpha_x \ge \lfloor \frac{n}{k}\rfloor$ for any vertex $x$.
\end{theorem}
\noindent {\bf Proof.} We only prove the sufficiency.
Suppose $(A,B)$ is a bipartition of $F=(V,E)$ with $|A|=a \ge
|B|=b$. Then $n=a+b$. Without loss of generality, we may assume
that $A$ has as few isolated vertices as possible. Let
$s_i=\lfloor \frac{n+i-1}{k} \rfloor$ for $1 \le i \le k$.  We
only need to partition $V$ into stable sets of size $s_1$, $s_2$,
$\ldots$, $s_k$, respectively. Choose the minimum index $j$ for
which $b \le \sum_{i=1}^j s_i$. If the inequality is an equality,
we can partition $V$ into desired stable sets. So, we now assume
that $\sum_{i=1}^{j-1}
 s_i < b < \sum_{i=1}^j s_i$.

{\bf Case 1.} $ 1 < j$.

Let $S$ be the set of $s=b - \sum_{i=1}^{j-1} s_i$ vertices of
lowest degrees in $B$. The number of edges between $S$ and $A$ is
then at most $s$ times the average degree of a vertex in $B$,
which is at most $\frac{n-1}{b}$. Therefore, $|N(S)| \le
\frac{s(n-1)}{b}<\frac{sn}{b}$ and then
$$
\textstyle{ |S \cup (A-N(S))| > s + a - \frac{sn}{b}
                            = \frac{(b-s)a}{b} \ge s_1, }
$$
since $b-s\ge s_1$ and $a\ge b$. Hence, $|S \cup (A-N(S))| \ge s_1
+ 1 \ge s_j$ and we can find a subset $S'$ of $A$ such that $S
\cup S'$ is a stable set of size $s_j$. In this case, the other
vertices can be properly partitioned to get an equitable
$k$-coloring of $F$.

{\bf Case 2.} $j=1$, i.e., $b < \lfloor \frac{n}{k} \rfloor$.

In this case, by the choice of $(A,B)$, we know that $A$ has no
isolated vertices. Denote $L$ the set of all leaves in $A$.  Then,
$|L| + 2 |A-L| \le \sum_{x\in A} \deg(x) \le n-1$ and so $|L| \ge
|L|+ |L|+2|A-L| -(n-1)=2a-n+1=a-b+1$.

We first choose a subset $S$ of $B$ such that the stable set
$(N(S)\cap L)\cup (B-S)$ has size at least
$\lceil\frac{n}{k}\rceil$.  Notice that since $k\ge 3$ and $b <
\lfloor \frac{n}{k} \rfloor$, we have $|L| \ge \lceil \frac{n}{k}
\rceil$. Hence, $B$ is such a candidate, while $\emptyset$ is not.
We may assume that $S$ is chosen so that $|S|$ is smallest. Choose
a vertex $v$ from $S$.  Then $|N(S-\{v\})\cap L|+ |B-(S-\{v\})|
\le \lceil\frac{n}{k} \rceil - 1$.

If the stable set $(N(B-S)\cap L) \cup S$ has size at least
$\lfloor \frac{n}{k}\rfloor$, then $A$ has two disjoint subsets
$S'$ and $S''$ such that $S'\cup(B-S)$ and $S''\cup S$ are two
stable sets of size $s_k$ and $s_1$, respectively. Hence the other
vertices can be properly partitioned to get an equitable
$k$-coloring of $F$. So, we may assume that $|N(B-S)\cap L|+|S|
\le \lfloor \frac{n}{k}\rfloor - 1$.
Adding the two inequalities gives $|L|-|N(v)\cap L|+b+1 \le
\lceil\frac{n}{k} \rceil + \lfloor \frac{n}{k}\rfloor -2$.
Consequently,
$$
\textstyle{  |N(v)\cap L| \ge |L|+b +1-\lceil\frac{n}{k} \rceil -
\lfloor \frac{n}{k}\rfloor +2 \ge a + 4 - \lceil\frac{n}{k} \rceil
- \lfloor \frac{n}{k}\rfloor.  }
$$

Since $\alpha_v \ge \lfloor\frac{n}{k}\rfloor$, there is a
$v$-stable set $R$ of size $\lfloor\frac{n}{k}\rfloor$.  We may
assume that $R$ is chosen so that $|R\cap B|$ is minimum.
If $R\cap B=\{v\}$, then as $|(N(v)\cap L)\cup (B-\{v\})|\ge a + 4
- \lceil\frac{n}{k} \rceil - \lfloor \frac{n}{k}\rfloor +b-1 \ge
\lceil\frac{n}{k} \rceil$ we can choose a subset $S'$ of $A$ such
that $S'\cup(B-\{v\})$ is a stable set of size
$\lceil\frac{n}{k}\rceil$. This and $R$ together with a proper
partition of other vertices give an equitable $k$-coloring of $F$.
Suppose $R\cap B$ has at least two vertices.  In this case, any
vertex $x\in L$ that is not in $R$ must be adjacent to some vertex
in $R\cap B$, for otherwise we can replace a vertex in $(R\cap B)
- \{v\}$ to get a $v$-stable set $R'$ of the size
$\lfloor\frac{n}{k}\rfloor$, but $|R'\cap B| < |R\cap B|$,
contradicting the choice of $R$.  Therefore, any vertex of $L$ is
either in $R$ or adjacent to some vertex in $R$.  Then $(B \cup L)
-R$ is a stable set of size at least
$b+(a-b+1)-\lfloor\frac{n}{k}\rfloor \ge \lceil\frac{n}{k}\rceil$.
Again, we are able to equitably $k$-color $F$. \qed

Notice that while it is easy to characterize equitable
$2$-colorability of a tree, it is slightly complicated for a
forest.  Suppose a forest $F$ of order $n$ has $r$ components,
each has order $n_i=a_i+b_i$ where $a_i$ and $b_i$ are the sizes
of its partite sets.  To check the equitable $2$-colorability of
$F$ is the same as to partition $\{1,2,\ldots, r\}$ into $I$ and
$J$ such that $\sum_{i\in I} a_i + \sum_{j\in J} b_j =
\lfloor\frac{n}{2}\rfloor$.

We close this note by raising the problem that how far can we go
from trees to chordal graphs on equitable colorability.

\noindent {\bf Acknowledgements.} The author thanks Professors
Chen and Lih for introducing the equitable colorability of graphs
to him.

\end{document}